\documentclass[11pt,leqno]{article}

\usepackage{amsmath,theorem,amssymb,amscd}

\long\def\comment#1\endcomment{}

\makeatletter
\begingroup
\gdef\th@dotted{\normalfont\itshape
  \def\@begintheorem##1##2{%
        \item[\hskip\labelsep \theorem@headerfont ##1\ ##2.]}%
\def\@opargbegintheorem##1##2##3{%
   \item[\hskip\labelsep \theorem@headerfont ##1\ ##2\ (##3).]}}
\endgroup
\makeatother

\theoremstyle{dotted}

\newtheorem{thm}{Theorem}[section]
\newtheorem{lemma}[thm]{Lemma}
\newtheorem{remark}[thm]{Remark}

\newtheorem{prop}[thm]{Proposition}

\makeatletter
\begingroup
\gdef\th@upshape{\normalfont
  \def\@begintheorem##1##2{%
        \item[\hskip\labelsep \theorem@headerfont ##1\ ##2.]}%
\def\@opargbegintheorem##1##2##3{%
   \item[\hskip\labelsep \theorem@headerfont ##1\ ##2\ (##3).]}}
\endgroup
\makeatother

\theoremstyle{upshape}

\newtheorem{defn}[thm]{Definition}
\newtheorem{rem}[thm]{Remark}

\makeatletter
\renewcommand{\subsection}{\@startsection{subsection}{2}{0pt}{-3ex
plus -1ex minus -0.2ex}{-2mm plus -0pt minus
-2pt}{\normalfont\bfseries}} 
\makeatother

\makeatletter
\@addtoreset{equation}{section}
\makeatother

\newcommand{\proof}[1][Proof.]{\smallskip\noindent{\em #1}}
\def\endproof{\hfill\ensuremath{\square}\par\medskip}

\def\eqref#1{\thetag{\ref{#1}}}

\renewcommand{\labelenumi}{{\normalfont(\roman{enumi})}}

\let\latexref=\ref
\def\ref#1{{\normalfont{\latexref{#1}}}}

\newcommand{\wt}{\widetilde}

\newcommand{\6}{\partial}

\newcommand{\cntrct}
{\hspace{2pt}\raisebox{1pt}{\text{$\lrcorner$}}\hspace{2pt}}

\newcommand{\hdot}{{\:\raisebox{3pt}{\text{\circle*{1.5}}}}}
\newcommand{\arrow}{{\:\longrightarrow\:}}
\newcommand{\C}{{\mathbb C}}
\newcommand{\R}{{\mathbb R}}

\newcommand{\LL}{{\cal L}}

\newcommand{\X}{{\mathfrak X}}

\newcommand{\calo}{{\cal O}}
\newcommand{\T}{{\cal T}}

\newcommand{\m}{{\mathfrak m}}

\newcommand{\Spec}{\operatorname{Spec}}
\newcommand{\Ob}{\operatorname{\sf Ob}}
\newcommand{\Sympl}{\operatorname{\sf Sympl}}
\newcommand{\Aut}{\operatorname{\sf Aut}}
\newcommand{\Diff}{\operatorname{\sf Diff}}
\newcommand{\Tr}{\operatorname{\sf Tr}}
\newcommand{\Hom}{\operatorname{Hom}}
\newcommand{\Ext}{\operatorname{Ext}}
\newcommand{\Ham}{\operatorname{{\cal H}\!{\it am}}}

\newcommand{\Per}{\operatorname{\sf Per}}
\newcommand{\perdom}{\operatorname{\sf Per}}
\newcommand{\cPer}{{\cal P}{\it er}}

\newcommand{\Def}{\operatorname{\sf Def}}
\newcommand{\cDef}{\operatorname{{\cal D}{\it ef}}}

\newcommand{\KS}{\operatorname{{\cal K}{\cal S}}}

\newcommand{\Ker}{\operatorname{Ker}}
\newcommand{\id}{\operatorname{\sf id}} 
\renewcommand{\dim}{\operatorname{\sf dim}} 

\title{Period map for non-compact holomorphically symplectic
manifolds} 

\author{D. Kaledin and M. Verbitsky}

\begin{document}

\maketitle

{\small 
\hspace{0.2\linewidth}
\begin{minipage}[t]{0.7\linewidth}
{\bf Abstract:} We study the deformations
of a holomorphic symplectic manifold $M$,
not necessarily compact, over a formal ring. We show 
(under some additional, but mild, assumptions on $M$) 
that the coarse deformation space exists and
is smooth, finite-dimensional and naturally
embedded into $H^2(M)$. For a holomorphic symplectic
manifold $M$ which satisfies
$H^1(\calo_M) = H^2(\calo_M)=0$, the
coarse moduli of formal deformations
is isomorphic  to $\Spec\C[[t_1, ..., t_n]]$,
where $t_1$, ... $t_n$ are coordinates in $H^2(M)$.
\end{minipage}
}

\tableofcontents


\section{Introduction}


Deformations and moduli of compact K\"ahler manifolds are a well studied
subject, dating back to Kodaira-Spencer \cite{_Kodaira_Spencer_}. 
The moduli of non-compact manifolds are rarely mentioned, mostly
because they are much harder to define and study. 

The work on the moduli of compact holomorphically symplectic
manifolds and Calabi-Yau manifolds is still far from the conclusion;
the local case is due to F.Bogomolov, A.Beauville, G.Tian,
A.Todorov, P.Deligne and Z.Ran (\cite{bogo}, \cite{_Beauville_},
\cite{tian}, \cite{todo}, \cite{ran}).  Extreme importance of this
subject is highlighted by thousands of papers on Mirror Symmetry,
which appeared since then.

In the non-compact case, some work in this direction was done by
M.Kon\-tse\-vich and S.Barannikov (\cite{_Barannikov_Kontse_}) and
others, but, for the most part, this territory is still
uncharted. However, there are many examples that suggest that at
least for some non-compact holomorphically symplectic manifolds a
good local deformation theory does exist. In particular, in the
well-studied case of smooth crepant resolutions of symplectic
quotient singularities in dim $2$ (the so-called Du Val points), a
likely candidate for the universal local deformation is provided by
the simultaneous resolution of Brieskorn \cite{brisk}.

In this paper we extend to the non-compact case the algebraic
construction of the local deformation space of a Calabi-Yau manifold
$M$, due to Z.Ran. Unfortunately, our results are valid only when
the manifold $M$ is holomorphically symplectic.

Our approach is essentially the same as the original approach of
Bogomolov. It is based on the so-called {\bf period map}. Instead of
deformations of a holomorphically symplectic manifold $M$, one
considers deformations of the pair $\langle M, \Omega\rangle$, where
$\Omega$ denotes the holomorphic symplectic form, that is, a nowhere
degenerate closed $(2,0)$-differential form.  Given a local
deformation $\pi:\; \wt M \arrow S$ of $\langle M, \Omega \rangle$
with a simply connected base, the cohomology of the individual
fibers of $\pi$ are identified by the Gauss-Manin connection. Taking
the cohomology class of the holomorphic symplectic form of each
fiber, one obtains a map $\Per:\; S \arrow H^2(M)$.  Bogomolov and
Beauville have shown that for $M$ compact, the map $\Per$ induces a
holomorphic immersion of the coarse marked moduli space $\cal M$ of
$\langle M, \Omega \rangle$ into $H^2(M)$. The image of $\Per$
belongs to a certain quadric ${\cal C}\subset H^2(M)$, cut by the
so-called {\bf Bogomolov-Beauville form}.  Moreover, the period map
$\Per:\; \cal M \arrow \cal C$ is locally an isomorphism.

Bogomolov extended these results to Calabi-Yau manifolds in an
unpublished I.H.E.S preprint (1982). In 1987, Tian and Todorov
published a different proof of Bogomolov's theorem. Their proof of
Bogomolov-Tian-Todorov theorem was based on Hodge theory. An
algebraic version of their arguments was proposed by Z.Ran
(\cite{ran}). Ran's argument uses the degeneration of the $E_2$-term
of the Dolbeault spectral sequence (proven in algebraic case by
Deligne and Illusie, \cite{_Deligne_Illusie_}).  However, this
spectral sequence is not degenerate in non-compact case, hence this
argument does not work for open Calabi-Yau manifolds.

We found that a version of Z.Ran's argument is valid for holomorphic
symplectic manifolds (under some additional, quite weak,
assumptions).  One can explain this heuristically as follows.  For a
complex manifold $M$, deformations are classified by $H^1(\T(M))$,
where $\T(M)$ is the tangent sheaf. When $M$ is Calabi-Yau, $\T(M)$
is isomorphic to $\Omega^{n-1}(M)$, where $n=\dim_\C M$. To show
that the deformations of $M$ have no obstructions, we would need to
prove that the $E_2$-term of the Dolbeault spectral sequence
degenerates in $H^1(\Omega^{n-1}(M))$. This is very far from truth
in the non-compact case.  However, if $M$ is holomorphic symplectic,
we have $\T(M) \cong \Omega^1(M)$, and instead of
$H^1(M,\Omega^{n-1}(M))$ we have to consider
$H^1(M,\Omega^1(M))$. The only differential in the spectral sequence
that maps into $H^1(M,\Omega^1(M))$ starts at $H^1(M,\calo_M)$.  If
we assume for simplicity that $H^i(M,\calo_M)=0$ for $i \geq 1$,
then this differential vanishes tautologically. The differentials
that {\em start} at $H^1(M,\Omega^1(M))$ can still be non-trivial,
but it turns out that they become irrelevant if instead of
deformations of $M$ one considers deformations of the pair $\langle
M,\Omega\rangle$.

Thus in the case $H^i(M,\calo_M) = 0$, $i \geq 1$ we obtain the
following result.

\begin{thm}
Let $M$ be a holomorphic symplectic manifold such that for every $i
\geq 1$ we have $H^i(\calo_M) = 0$.  Then there exists a coarse
moduli space of formal deformations
$$
\pi:\; \wt{M} \arrow Spl(M, \Omega)
$$
of pairs $\langle M,\Omega \rangle$. Moreoved, the period map
$\Per:\; Spl(M, \Omega)\arrow H^2(M)$ gives an isomorphism of
$Spl(M, \Omega)$ and the formal completion of $H^2(M)$ in
$[\Omega]\in H^2(M)$.\endproof
\end{thm}

There is a more general version of this result (Theorem \ref{main}),
which works for a larger class of non-compact holomorphically
symplectic manifolds. Here we also obtain a coarse moduli space
$Spl(M,\Omega)$, which is smooth and finite-dimensional, but it is
no longer isomorphic to $H^2(M)$. However, the period map remains an
immersion.

\bigskip

We will now give a semi-rigorous sketch of the proof of this
result. First, consider an easy but important example of affine
holomorphic symplectic manifold $M$.  We have $H^1(\T(M))=0$, hence any
first-order complex deformation of $M$ is trifial. Using induction,
it is easy to show that any formal complex deformation of $M$ is
also trivial, that is, for any formal deformation $\pi:\; \wt{M}
\arrow S$ of $M$, we have $\wt{M} \cong S\times M$.  However, a
symplectic deformation needs not to be trivial, because we may have
non-trivial variations of holomorphic symplectic structure. A formal
deformation of the pair $\langle M, \Omega\rangle $ is determined by
the deformation of a closed $(2,0)$-form $\Omega$.

By Grothendieck's theorem, the topological cohomology of $M$ is
isomorphic to the hypercohomology of the algebraic de Rham complex
\begin{equation}\label{_de_Rham_Equation_} 
\begin{CD}
0 @>>> \calo_M @>\6>> \Omega^1 M
@>\6>> \Omega^2 M @>\6>> \dots
\end{CD}
\end{equation}
Since $M$ is affine, $H^i(\Omega^j (M))=0$ for $i>0$.  Therefore,
$H^2(M)$ is isomorphic to the second cohomology of the complex
\eqref{_de_Rham_Equation_}.  A first order infinitesimal
automorphism of $M$ is given by the section of $TM$, which is
isomorphic to $\Omega^1(M)$. A (co-)vector field $\gamma\in
\Omega^1(M) \cong \T(M)$ acts on $\Omega^i M$ by the Lie derivative.
Therefore, $\gamma$ acts on a closed form $\Omega\in \Omega^2(M)$ by
adding $d\gamma$, and the first-order deformations of $\langle M,
\Omega \rangle$ are classified by closed 2-forms up to exact
2-forms.  Using induction, it is easy to check that this is true in
any order. Therefore, the period map gives an isomorphism of the
coarse moduli of $\langle M, \Omega \rangle$ and the formal
neighbourhood of the class $[\Omega]$ in $H^2(M)$.

The general proof is deduced from the affine version as follows.

Let $S$ be a spectrum of an Artin ring over $\C$.  We consider the
deformations of $\langle M, \Omega \rangle$ over $S$ as a stack of
groupoids over $M$. We introduce another stack, called {\bf
Kodaira-Spencer stack}, which -- whenever it is defined --
classifies the maps from $M$ to $H^2(M)$. Using the affine version
of the main theorem, we show that these stacks are equivalent over
affine open subsets of $M$. This immediately implies that these
stacks are equivalent over the whole $M$, which concludes the study
of the formal deformations of $\langle M, \Omega\rangle$ and shows
that these deformations are classified by the maps form the base of
deformation to $H^2(M)$.

The above account is greatly simplified. To make it at least
approximately workable, we have to do big technical adjustments. The
problems are twofold: first of all, the Kodaira-Spencer groupoid is
not defined in a general situation -- we have to go step-by-step
through what we call {\em elementary} base extensions to define it
properly. Secondly, the Kodaira-Spencer groupoid classifies not the
maps from the base $S$ to $H^2(M)$, but rather a certain maps of
complexes of sheaves, which are reduced to the maps from $S$ to
$H^2(M)$ when $H^1(\calo_M)= H^2(\calo_M)=0$.

Now we give a more precise version of the definition of the
Kodaira-Spencer stack. Let $S$ be an Artin scheme over $\C$, and
$S_0\subset S$ a closed subscheme defined by an ideal $I\subset
\calo_S$, such that $I^2=0$. Assume that the ideal $I$ is
sufficiently small (see Definition~\ref{elem.ext} for the precise
condition on $I$).  Fix a deformation $\pi_0:\; \wt M_{S_0}
\arrow S_0$ of $\langle M, \Omega \rangle$.  Consider the set
$\Def(\wt M_{S_0}, S)$ of all deformations $\wt M_S$ of
$\langle M, \Omega \rangle$ over $S$, equipped with an isomorphism
\[ 
\wt M_S \times_{S} S_0 \arrow \wt M_{S_0}.
\]
Clearly,  $\Def(\wt M_{S_0}, S)$ is a stack of groupoids over
$M$. Consider a truncated relative 
de Rham complex $F^1 \Omega^*(\wt M_{S_0}/S_0)$,
\[ 
\begin{CD}
  \Omega^1 (\wt M_{S_0}/S_0) @>\6>>
  \Omega^2 (\wt M_{S_0}/S_0) @>\6>> \dots
\end{CD}
\]
A holomorphic symplectic structure on the deformation $\wt M_{S_0}$
defines a morphism of complexes 
$\pi_0^* \calo_{S_0}[2] \arrow F^1 \Omega^*(\wt M_{S_0}/S_0)$.\
Taking its derivative along the Gauss-Manin connection on 
$S_0$, we obtain the so-called
Kodaira-Spencer map of the deformation $\wt M_{S_0}$
\[ 
 \pi_0^* TS_0[2] \stackrel {\theta_0}\arrow F^1 \Omega^*(\wt
 M_{S_0}/S_0).
\]
The Kodaira-Spencer stack $KS(\wt M_{S_0}, S)$ is defined as
follows. The {\em objects} of $KS(\wt{M}_{S_0},S)$ are all morphisms
of complexes
\[ 
\pi_0^* TS\otimes {\calo_S} \calo_{S_0}[2] 
    \stackrel \theta \arrow F^1 \Omega^*(\wt M_{S_0}/S_0)
\]
such that the restriction of $\theta$ to $TS_0$ is equal to
$\theta_0$. The {\em morphisms} between any two such $\theta_1$,
$\theta_2$ are chain homotopies between them -- that is, maps 
\[ 
\pi_0^* TS\otimes {\calo_S} \calo_{S_0}[1] 
    \stackrel \gamma \arrow F^1 \Omega^*(\wt M_{S_0}/S_0)
\]
satisfying $d\gamma = \theta_1 - \theta_2$. This obviously defines a
groupoid. The definition of the period map from $\Def (\wt
M_{S_0}, S)$ to $KS(\wt M_{S_0}, S)$ is straightforward; using
the deformation theory for affine $M$, we show that this is an
isomorphism.  In other words, there are exactly as many ways to
extend the deformation from $S_0$ to $S$ as there are ways to extend
the corresponding morphism to $H^2(M)$. Using induction, this leads
immediately to the classification result for the deformations stated
above.

\bigskip

The paper is organized as follows. We start with a semi-heuristic
analytic proof of the main theorem obtained by the use of the
Cartan-Maurer equation, as suggested by Kontsevich and
Barannikov. This proof cannot be made precise in the non-compact
situation, but we have decided to include it anyway, since it
exhibits nicely the general idea behind the rigorous proof. Then we
turn to purely algebaric methods. Section 3 contains the necessary
definitions and the precise statement of our result in the strongest
form. In Section 4, we introduce the approriate symplectic version
of the Kodaira-Spencer class and define elementary extensions. In
Section 5, we rework our construction in the language of stacks and
prove the main theorem. Finally, Section 6 is a short postface which
explains how our results are related to the known facts (in
particular, to the results of Z. Ran).


\section{Deformations in mixed formal complex-analytic category}
\label{_Koncevich_Section_}


One of the way to study deformations us through the $DG$-algebra
approach suggested by Kontsevich and Barannikov
(\cite{_Barannikov_}, \cite{_Barannikov_Kontse_}). The main
advantage of this method is that it foregoes the tedious
step-by-step constructions of Grothendieck's local deformations
theory, and gives the results in the form of an explicit power
series. Unfortunately, it is not quite as useful in non-compact
case, when we have no means to check that these series converge.

The deformations one obtains from the Dolbeault complex lie a weird
mixed formal-complex-analytic category. To obtain something
definite, one needs some kind of integrability-type
conditions. Kontsevich and Barannikov work with compact manifolds,
so that in their situation this is not a problem -- integrability
can be obtained directly from functional analysis. However, since we
study open manifolds, the functional analysis does not help, and we
stay in the mostly useless mixed category.

Nevertheless, Kontsevich's approach (which goes back to Kodaira) is
beautiful and quite useful as a heuristic tool. Therefore we decided
to express some of our results in this language before going to the
rigourous step-by-step proof.

Since the $DG$-algebra approach is used only for
heuristics, this section will be quite sketchy;
a cursory knowledge of \cite{_Barannikov_Kontse_}
is required.

We start with a review of the deformation theory as it is given in
\cite{_Barannikov_} and \cite{_Barannikov_Kontse_}.
Let $M$ be a complex manifold,
$M_\R$ the underlying real analytic manifold 
and $M_\R[[t]]:= M_\R\times \Spec(\C[[t]])$
the ``mixed formal-real analytic'' manifold obtained
as a product of $M_\R$ and the formal disk
$\Delta:= \Spec(\C[[t]])$. 
Using $DG$-algebras, one may classify the complex
deformations of $M$ over $\Delta$, that is, the
complex structures $J$ on $M_\R[[t]]$ 
such that the zero fiber of $(M_\R[[t]], J)$ is isomorphic
to $M$. It is well known that such deformations are classified by
the solutions of the Maurer-Cartan equation
\begin{equation}\label{_Maurer_Cartan_Equation_} 
  \bar \6 \gamma(t) = - \frac{1}{2}[ \gamma(t), \gamma(t)],
\end{equation}
where $\gamma(t)\in \Lambda^{0,1}(TM)[[t]]$
is a $C[[t]]$-valued $(0,1)$-form with coefficients
in holomorphic vector fields, and
\[ [\cdot,\cdot] :\; \Lambda^{0,1}(TM)[[t]]\times \Lambda^{0,1}(TM)[[t]]
   \arrow \Lambda^{0,2}(TM)[[t]]
\]
is the Schouten bracket.

\medskip

We explain in a few words how the complex structures
are classified by the solutions of Maurer-Cartan.

Consider the sheaf $A^{\hdot,\hdot}=\Lambda^{\hdot,\hdot}(M)$ as an abstract
sheaf of algebras. A complex structure on $M$ is defined
by an identification between $\Lambda^1(M)$ and 
$A^{0,1}$. Since $\Lambda^1(M)$ classifies the
derivations of the sheaf of smooth functions,
to give an identification $\Lambda^1(M)\cong A^{0,1}$
is the same as to give a derivation
$\bar\6:\; C^\infty(M) \arrow A^{0,1}$.
The difference between two such operators 
is given by $\gamma\in \Lambda^{0,1}(TM)$.
An integrability condition is written as 
$(\bar\6+\gamma)^2=0$, which is rewritten as
the Maurer-Cartan equation. Deformations
are equivalent if the corresponding operators
are exchanged  by an automorphism of $A^{\hdot,\hdot}$.

We are going to write a similar interpretation for
holomorphic symplectic deformations.

\medskip

Fix a holomorphic symplectic form $\Omega\in A^{2,0}$.
A holomorphic symplectic deformation
is defined by an operator
\[ \bar\6+\gamma:\;  C^\infty(M) \arrow A^{0,1},\ \ \ 
   \gamma(\Omega)=0
\]
such that $(\bar\6+\gamma)^2$, or, equivalently,
$\gamma$ satisfies the Maurer-Cartan 
equation \eqref{_Maurer_Cartan_Equation_}.
Deformations are equivalent if the corresponding operators
are exchanged  by an automorphism of $A^{\hdot,\hdot}$
preserving $\Omega$.

\medskip

In this spirit, Barannikov 
\cite{_Barannikov_} 
describes the
deformations of Calabi-Yau manifolds, 
with $\Omega$ a nowhere degenerate section of
the sheaf of top-degree $(p,0)$-forms.

\medskip

For the holomorphic symplectic manifolds, the deformations are
described in terms the sheaf $\Ham$ of 
holomorphic Hamiltonian vector fields.
First of all, automorphisms of $A^{\hdot,\hdot}$
preserving $\Omega$ correspond to 
the Hamiltonian vector fields $\beta \in \Ham(M)$.
Secondly, the condition $\gamma\Omega=0$ means that
$\gamma\in \Lambda^{0,1}(\Ham(M))$.

We obtain that the first-order deformations of
$M$ are described by the cohomology of the sheaf 
of Hamiltonian vector fields. To prove that
the deformations of $M$ are unobstructed
(that is, every one-parametric first-order deformation is extended
to a deformation over a formal disk), we need to show
the following. Let $\gamma_1\in \Lambda^{0,1}(\Ham(M))$
be a $\bar\6$-closed $(0,1)$-form with coefficients
in $\Ham$. Then there exists a series
$\gamma_2, \gamma_3 , \dots \in \Lambda^{0,1}(\Ham(M))$
such that
\[ (\bar \6 +t \gamma_1 + t^2 \gamma_2 + \dots)^2 =0,
\]
that is, the formal series 
$\gamma_1 +t\gamma_2 + t^2 \gamma_3+\dots$
satisfies the Maurer-Cartan 
equation \eqref{_Maurer_Cartan_Equation_}.
In this setting, Maurer-Cartan can be written as follows:
\begin{equation}\label{_Maurer_Cartan_serie_Equation_}
   \bar\6\gamma_{n+1}=-\frac{1}{2}\sum_{i+j=n}[\gamma_i, \gamma_j],
   \ \ \  n=2, 3, \dots 
\end{equation}
Consider the following commutative diagram
\begin{equation}\label{_polotno_Equation_}
\begin{CD}
@>{\bar\6}>> \Lambda^{2,1} @>{\bar\6}>> \Lambda^{2,2}@>{\bar\6}>>\\
& & @AA\6 A @AA\6 A\\[3mm]
@>{\bar\6}>>  \Lambda^{1,1} @>{\bar\6}>> \Lambda^{1,2}@>{\bar\6}>>\\
&& @AA\6 A @AA\6 A\\[3mm]
@>{\bar\6}>>  \Lambda^{0,1} @>{\bar\6}>> \Lambda^{0,2}@>{\bar\6}>>
\end{CD}
\end{equation}
Identifying $TM$ and $\Lambda^{1,0}$, we can 
realize the polyvector fields as differential $(p,q)$-forms.
In particular, the solutions $\gamma_i$ of 
\eqref{_Maurer_Cartan_serie_Equation_}
are considered now as sections of $\Lambda^{1,1}$. 
The Hamiltonian condition is written
as $\6 \gamma_i =0$, where
$\6:\; \Lambda^{1,1}\arrow \Lambda^{2,1}$
is an operator in \eqref{_polotno_Equation_}.
Consider the 
``symplectic Hodge operator''
\[ \Lambda:\; \Lambda^{p,q} \arrow \Lambda^{p-2,q} 
\]
which is adjoint to the multiplication by the
holomorphic symplectic form $\Omega$ via
the non-degenerate product defined by $\Omega$.
The Schouten bracket is written in terms of
$\Lambda$ and $\6$ as follows.

\begin{lemma}\mbox{}\quad{\bf (Tian-Todorov lemma for holomorphic
symplectic manifolds.)} 
\label{_Tian_Todo_Lemma_}
Let $\gamma, \gamma'\in \Lambda^{1,1}$,
and let $[\gamma, \gamma']\in \Lambda^{2,1}$
denote the Schouten bracket (we identify
$TM$ and $\Lambda^{1,0}$ using the holomorphic
symplectic form). Then
\[ [\gamma, \gamma'] =
   \6\Lambda(\gamma\wedge \gamma')
   - \Lambda(\6\gamma\wedge \gamma') -\Lambda(\gamma\wedge \6\gamma').
\]
\end{lemma}
\proof{} 
This is a symplectic version of the standard
Tian-Todorov lemma, and it is proven in exactly the same fashion
as the usual one (\cite{todo}). 
\endproof

The main result of this Section is the following theorem,
which states that the deformations of holomorphic symplectic
manifold are unobstructed, given that
$H^{2}(\calo_M)=0$. This assumptions is
not essential, but it makes the statements much simpler.

\begin{thm}\label{defo_through_DG_Theorem_}
Let $M$ be a complex holomorphic symplectic
manifold, and let $\gamma_1\in \Lambda^{0,1}(\Ham)$ be
a $\bar\6$-closed 1-form with coefficients in Hamiltonian vector fields.
Assume that $H^{2}(\calo_M)=0$, that is, the 
holomorphic cohomology of the structure sheaf vanish.
Then \eqref{_Maurer_Cartan_serie_Equation_}
has a solution.
\end{thm}

\proof{} 
Let us write \eqref{_Maurer_Cartan_serie_Equation_}
in terms of the diagram \eqref{_polotno_Equation_},
using the isomorphism $TM\cong \Lambda^{1,0}$.
Using induction, we may assume that
\[ \bar\6 \sum_{i+j=n}[ \gamma_i, \gamma_j] =
   \sum_{i+j+k=n}[ \gamma_i, [\gamma_j,\gamma_k]] +
   [[ \gamma_i,\gamma_j,]\gamma_k],
\]
which is equal zero by Jacoby identity.
Therefore, 
\begin{equation}\label{_sum_[gamma,]_Equation_}
  \sum_{i+j=n}[ \gamma_i, \gamma_j]
\end{equation}
is $\bar \6$-closed; to solve 
\eqref{_Maurer_Cartan_serie_Equation_}
we need to show that it is $\bar\6$-exact,
that is, to prove that it represents a zero
class in the cohomology $H^2(\Ham(M))$.
The Hamiltonian vector fields are identified
with $\6$-closed $(1,0)$-forms. 
Poincare lemma gives an exact sequence
of sheaves
\begin{equation}\label{_Hami_exact_seque_Equation_}
   0 \arrow \C \arrow \calo_M \stackrel \6 \arrow \Ham(M) \arrow 0,
\end{equation}
Consider the piece
\[ \arrow H^2(\calo_M)\arrow H^2(\Ham) \stackrel h \arrow H^3(M) \arrow 
\]
of the corresponding long exact sequence.
By Lemma \ref{_Tian_Todo_Lemma_}, the
expression \eqref{_sum_[gamma,]_Equation_} is $\6$-exact;
therefore, it represents zero in $H^3(M)$. 
By our assumptions, $H^2(\calo_M)=0$, and therefore,
the map $h:\;  H^2(\Ham) \arrow H^3(M)$
is an embedding. This proves that the sum
\eqref{_sum_[gamma,]_Equation_} represents zero
in $H^2(\Ham(M))$, and there exists
$\gamma_{n}\in \Lambda^{0,1}(\Ham(M))$
such that $\bar\6\gamma_n = \sum_{i+j=n}[ \gamma_i, \gamma_j]$.
We proved Theorem \ref{defo_through_DG_Theorem_}.
\endproof

In Theorem \ref{defo_through_DG_Theorem_}
we identified the deformation space of 
$\langle M, \Omega \rangle$ with the cohomology of $\Ham(M)$.
These cohomology spaces are very easy to write down explicitly.
Writing the long exact sequence
corresponding to \eqref{_Hami_exact_seque_Equation_}, we obtain
\[ H^1(M) \arrow H^1(\calo_M)\arrow H^1(\Ham(M)) \arrow H^2(M)\arrow  
   H^2(\calo_M).
\]
This gives an identification of the deformation space with
$H^2(M)$ when $H^1(\calo_M)=H^2(\calo_M)=0$,
and, more generally, allows one to express 
the holomorphic symplectic deformations through
$H^i(M)$, $H^i(\calo_M)$ $(i=1,2)$.

\bigskip

In the remaining part of the paper, we combine the intuition
of the $DG$-algebra approach with the hard science
of Grothendieck's local deformations theory, and
obtain essentially the same results in a much
more rigorous setting.

\section{Statement of the results.} 

\subsection{Admissible manifolds.} To begin with, we will describe the
restrictions that we need to impose on the given holomorphically
symplectic manifold $X$.

Let $X$ be a smooth algebraic manifold over $\C$. By a {\em
deformation} of $X$ over a pointed scheme $\langle S, o \in S
\rangle$ we will understand a scheme $\wt{X}/S$ smooth over $S$ and
equipped with an isomorphism $o \times_S \wt{X} \cong X$. We will
only consider deformations over spectra $S = \Spec A$ of local Artin
$\C$-algebra $A$, so that the fixed point is given by the maximal
ideal $\m \subset A$. Unless otherwise mentioned, all deformations
will be assumed to be of this type.

Recall that the de Rham complex $\Omega^\hdot(X)$ is equipped with
the Hodge, a.k.a. stupid filtration $F^\hdot\Omega^\hdot(X)$ given by
$$
F^i\Omega^j(X) = \begin{cases} \Omega^j(X), &\quad j \geq i,\\
0, &\quad \text{otherwise}. \end{cases}
$$
This filtration is also defined, and by the same formula, for the
relative de Rham complex $\Omega^\hdot(\wt{X}/S)$ of an arbitrary
deformation $\wt{X}/S$.

For symplectic deformations, it is the first term
$F^1\Omega^\hdot(X)$ of the Hodge filtration that plays the crucial
role. 

The complex $F^1\Omega^\hdot(X)$ can be included in an obvious exact
triangle
\begin{equation}\label{tria}
\begin{CD}
F^1\Omega^\hdot(X) @>>> \Omega^\hdot(X) @>>> \calo_X @>>> 
\end{CD}
\end{equation}
where $\calo_X$ is the structure sheaf of the manifold $X$. This
triangle induces an exact triangle on cohomology. 

\begin{defn}\label{adm}
A smooth algebraic manifold $X$ over $\C$ is called {\em admissible}
if for any deformation $\pi:\wt{X} \to S$ the relative cohomology sheaf
$$
R^2\pi_*F^1\Omega^\hdot(\wt{X}/S)
$$
is a flat sheaf on $S$ and the canonical map
$$
R^2\pi_*F^1\Omega^\hdot(\wt{X}/S) \to R^2\pi_*\Omega^\hdot(\wt{X}/S)
$$
is injective. 
\end{defn}

This definition is pretty technical, because it is given in the most
general form. For all practical applications that we see at the
moment, it suffices to assume the stronger condition on $X$ provided
by the following easy lemma.

\begin{lemma}\label{cmp}
Let $X$ be a smooth complex algebraic manifold. If for all $p \geq
1$ the canonical map
$$
H^p(X,\C) \to H^p(X,\calo_X)
$$
is surjective, then the manifold $X$ is admissible.
\end{lemma}

\proof{} Let $\pi:\wt{X} \to S$ be an arbitrary deformation. The
existence of the Gauss-Manin connection implies that
$$
R^p\pi_*\Omega^\hdot(\wt{X}/S) \cong H^p(X,\C) \otimes \calo_S
$$
for every $p \geq 0$. We will prove that 
\begin{enumerate}
\renewcommand{\labelenumi}{(A)}
\item for every $p \geq 2$, the canonical map
$R^p\pi_*\Omega^\hdot(\wt{X}/S) \to R^p\pi_*\calo(\wt{X})$ is
surjective and the sheaf $R^p\pi_*\calo(\wt{X})$ is flat on $S$.
\end{enumerate}
Use downward induction on $p$, starting with an arbitrary $p > 2\dim
X$. Assume (A) proved for all $p > k$. Denote by $i:o \to S$ the
embedding of the base point. Consider the restriction of the
cohomology exact triangle induced by \eqref{tria} to the base point
$o \in S$. Then base change and the Nakayma Lemma immediately imply
that the map
$$
R^k\pi_*\Omega^\hdot(\wt{X}/S) \to R^k\pi_*\calo(\wt{X})
$$
is surjective. This together with the inductive assumption implies
that
\begin{enumerate}
\renewcommand{\labelenumi}{(B)}
\item the canonical map $R^p\pi_*F^1\Omega^\hdot(\wt{X}/S) \to
R^p\pi_*\Omega^\hdot(\wt{X}/S)$ is injective and the sheaf
$R^p\pi_*F^1\Omega^\hdot(\wt{X}/S)$ if flat over $S$ for every $p
\geq k+1$.
\end{enumerate}
Therefore if the sheaf $R^k\pi_*\calo(\wt{X})$ is not flat, then the
coboundary map
$$
L^1i^*R^k\pi_*\calo(\wt{X})) \to
i^*R^k\pi_*F^1\Omega^\hdot(\wt{X}/S)
$$
is not zero, which contradicts the assumption. Carrying the
induction half-step further, we derive that the map
$$
R^1\pi_*\Omega^\hdot(\wt{X}/S) \to R^1\pi_*\calo(\wt{X})
$$
is surjective. Therefore (B) also holds for $p=2$, which proves the
lemma.  
\endproof

This lemma shows that a manifold $X$ is admissible in the following
cases: 
\begin{enumerate}\label{lst}
\item $X$ is compact (Hodge theory).
\item $X$ is affine ($H^i(X,\calo_X) = 0$ for $i \geq 1$).
\item More generally, $X$ admits a proper generically one-to-one map
$\pi:X \to Y$ into an affine variety $Y$ with rational singularities
(again $H^i(X,\calo_X) = 0$ for $i \geq 1$).
\item $X$ admits a proper generically one-to-one map $\pi:X \to Y$
into an affine variety $Y$ and has trivial canonical bundle $K_X$
($H^i(X,\calo_X) = H^i(X,K_X) = 0$ for $i > 0$. \footnote{This
follows immediately from the
Grauert-Riemenschneider Vanishing Theorem, \cite{_Grauert_}.}
\end{enumerate}

If $X$ is an admissible manifold, one can choose a splitting
$H^2(X,\C) \to H^2(X,F^1\Omega^\hdot(X))$ of the canonical
surjection $H^2(X,F^1\Omega^\hdot(X)) \to H^2(X,\C)$. Together with
the Gauss-Manin connection, this splitting defines an isomorphism
$$
R^2\pi_*F^1\Omega^\hdot(\wt{X}/S) \cong H^2(X,F^1\Omega^\hdot(X))
\otimes \calo_S
$$
for every deformation $\pi:\wt{X} \to S$. We will always assume
given such a splitting, keeping in mind that this introduces into
our constructions an element of choice. Note that this choice does
{\em not} appear in the case of the affine $X$, since in this case
for every deformation we have $R^2\pi_*\calo(\wt{X}) = 0$ and
$R^2\pi_*F^1\Omega^\hdot(\wt{X}/S) \cong
R^2\pi_*\Omega^\hdot(\wt{X}/S)$.

\subsection{Symplectic deformations and the period map.}

Let $X$ be an admissible manifold. Assume from now on that the
manifold $X$ is equipped with a non-degenerate closed $2$-form
$\Omega \in \Omega^2(X)$.

\begin{defn}\label{def.defn}
A {\em symplectic deformation} $\wt{X}/S$ of the symplectic manifold
$X$ over a base $S$ is a usual deformation $\pi:\wt{X} \to S$
equipped with a closed relative $2$-form $\Omega \in
\Omega^2(\wt{X}/S)$ which becomes the given $2$-form under the
isomorphism $o \times_S \wt{X} \cong X$.
\end{defn}

For every local Artin scheme $S = \Spec A$, we will denote by
$\Def(X,S)$ or simply by $\Def(S)$ the set of isomorphism classes of
symplectic deformations $\wt{X}/S$ of $X$ over $S$.

Choose once and for all a splitting $H^2(X,\C) \to
H^2(X,F^1\Omega^\hdot(X))$, so that for every deformation
$\pi:\wt{X} \to S$ we have an isomorphism
$$
R^2F^1\Omega^\hdot(\wt{X}/S) \cong H^2(X,F^1\Omega^\hdot(X)) \otimes
\calo_S.
$$
If the deformation $\wt{X}/S$ is symplectic, then the relative
$2$-form $\Omega \in \Omega^2(\wt{X}/S)$ defines a canonical
cohomology class
$$
[\Omega] \in H^2(X,F^1\Omega^\hdot(X)) \otimes \calo_S.
$$
This class gives a scheme map
$$
\Per(\wt{X}):S \to Tot(H^2(X,F^1\Omega^\hdot(X)),
$$
where $Tot(H^2(X,F^1\Omega^\hdot(X))$
denotes the total space of $H^2(X,F^1\Omega^\hdot(X))$
considered as a scheme. Further on, we shall simplify
notation by omitting ``$Tot$''.

\begin{defn}\label{perdom.defn}
The {\em period domain} of the admissible symplectic manifold $X$ is
the completion of the vector space $H^2(X,F^1\Omega^\hdot(X))$ near
the point $[\Omega] \in H^2(X,F^1\Omega^\hdot(X))$ corresponding to
the symplectic form $\Omega \in \Omega^2(X)$.

The map $\Per(\wt{X})$ is called {\em the period map} of the deformation
$\wt{X}/S$. 
\end{defn}

Note that the period map by definition maps $S$ into the period
domain $\perdom \subset H^2(X,F^1\Omega^\hdot(X))$.

\begin{rem}
This definition of the period domain is essentially cheating: in
Bogomolov's theory the period domain is not a formal scheme, but a
globally (and non-trivially) defined quadric in the projectivization
of the vector space $H^2(X,\C)$. However, since we work only with
infinitesemal deformation, Definition~\ref{perdom.defn} is
sufficient for our purposes.
\end{rem}

For any local Artin scheme $S$, taking the period map defines a map
$$
\Per:\Def(S) \to \perdom(S)
$$
from the set of deformation classes over $S$ to the set of
$S$-points of the formal scheme $\perdom$. This map (which we will,
by abuse of the language, also call the
period map) is functorial in $S$ (where $\Def(S)$ is considered as a
functor by taking pullbacks). We can now state our main result.

\begin{thm}\label{main}
Let $X$ be an admissible manifold equipped with a symplectic
$2$-form $\Omega \in \Omega^2(X)$. Then for any local Artin scheme
$S$, the period map
$$
\Per:\Def(S) \to \perdom(S)
$$
induces a set bijection
from the set of isomorphism classes of symplectic deformations
$\wt{X}/S$ to the set of $S$-points of the period domain $\perdom$. 

In particular, there exists a (formal) symplectic deformation
$\X/\perdom$ such that any deformation $\wt{X}/S$ is isomorphic to
the pullback of $\X$ by means of the period map $\Per(\wt{X}):S \to
\perdom$.
\end{thm}

\begin{remark}
Let $M$ satisfy $H^i(M,\calo_N)=0$ for $i \geq 1$.  Then Theorem
\ref{main} holds. Moreover, the period map $\Per:\Def(S) \to H^2(M)$
is locally an isomorphism.
\end{remark}

\begin{remark} For $M$ affine, we obviously have $H^i(M,\calo_M)=0$
for $i \geq 1$. $H^1(\calo_N)= H^2(\calo_M)=0$ obviously holds.  In
this case, Theorem \ref{main} can be proved by a standard inductive
argument (see Introduction).
\end{remark}

\section{Elementary extensions and the\! Kodaira-Spencer class.}

\subsection{Elementary extensions.}

To prove Theorem~\ref{main}, we will use induction on the length of
the local Artin algebra $A = \calo(S)$. To set up the induction, we
introduce the following.

\begin{defn}\label{rel.def.defn}
Let $S_0 \subset S$ be a closed embedding of local Artin schemes,
and let $\wt{X}_0/S_0$ be a symplectic deformation of a
holomorphically symplectic manifold $X$ over $S_0$. 

Then by $\Def(\wt{X}_0,S)$ we will denote the set of isomorphism
classes of symplectic deformations $\wt{X}/S$ equipped with a
symplectic isomorphism $\wt{X} \otimes_S S_0 \cong \wt{X}_0$.
\end{defn}

This is consistent with our earlier notation. Indeed, by
Definition~\ref{def.defn}, every symplectic deformation is
canonically trivialized over the base point $o \subset S$. Therefore
$\Def(X,S)$ in the sense of Definition~\ref{rel.def.defn} is still
the set of isomorphism classes of all symplectic deformations of the
manifold $X$.

The corresponding notion on the ``period'' side is the following.

\begin{defn}
Let $p_0:S_0 \to \perdom$ be a map from the closed subscheme $S_0
\subset S$ to the period domain $\perdom$. Then by $\Per(p_0,S)$ we
will denote the set of all maps $f:S \to \perdom$ such that
$F|_{S_0} = p_0$.
\end{defn}

It is customary in deformation theory to prove theorems
step-by-step, starting with the case of square-zero
extensions. However, we will need a slightly smaller class of
extensions $S_0 \subset S$. Namely, let $A$ be a local Artin
algebra, and let $I \subset A$ be a square-zero ideal, so that $I^2
= 0$. Then we have the usual exact sequence of the modules of
K\"ahler differentials over $\C$
\begin{equation}\label{diff}
\begin{CD}
I @>>> \Omega^1(A)/I @>>> \Omega^1(A/I) @>>> 0.
\end{CD}
\end{equation}
\begin{defn}\label{elem.ext}
The extension $\Spec A/I \subset \Spec A$ will be called {\em
elementary} if the sequence \eqref{diff} is also exact on the left.
\end{defn}

The following lemma immediately implies that every local Artin
scheme $S$ admits a filtration $o \subset S_0 \subset \ldots \subset
S_k = S$ such that all extensions $S_i \subset S_{i+1}$ are
elementary. 

\begin{lemma}
Let $A$ be a local Artin algebra with the maximal ideal $\m$. Assume
that $\m^{p+1} = 0$, while $\m^p \neq 0$. Then the extension $\Spec
A/\m^p \subset \Spec A$ is elementary.
\end{lemma}

\proof{} We have to prove that the canonical map $\m^p \to
\Omega^1(A)/\m^p$ is injective. It suffices to prove this for $A$
replaced with its associated graded quotient with respect to the
$\m$-adic filtration. Thus we can assume that $A$ and $\Omega^1(A)$
are graded. The map $\m^p \to \Omega^1(A)/\m^p$ is the composition
of the de Rham differential $d:\m^p \to \Omega^1(A)$ and the
projection $\Omega^1(A) \to \Omega^1(A)/\m^p$. Since $d$ is
obviously injective, it suffices to prove that $d(\m^p) \cap \m^p
\cdot \Omega^1(A) = 0$. But this is trivial: $d(\m^p)$ has degree
$p$ with respect to the grading on $\Omega^1(A)$, while $\m^p \cdot
\Omega^1(A)$ is of degree $(p+1)$.
\endproof

This lemma reduces Theorem~\ref{main} to the following claim.

\begin{prop}\label{rel.main}
Let $X$ be an admissible symplectic manifold.  Let $S_0 \subset S$
be an elementary extension of local Artin schemes, and let
$\wt{X}_0/S_0$ be an arbitrary symplectic deformation of the
manifold $X$. Denote by $p_0:S_0 \to \perdom$ the period map of the
deformation $\wt{X}_0/S_0$.

Then the period map
$$
\Per:\Def(\wt{X}_0,S) \to \perdom(p_0,S)
$$
is as isomorphism.
\end{prop}

\subsection{The Kodaira-Spencer class.} 

In order to start the proof of Proposition~\ref{rel.main}, we will
need a convenient description of the set $\Per(p_0,S)$. To give such
a description, we will consider not the period map itself, but its
differential.

\begin{defn}\label{ks.class.defn}
Let $\wt{X}/S$ be a symplectic deformation of a symplectic manifold
$X$. Then the class
$$
\theta = \nabla\Omega \in H^2(\wt{X},F^1\Omega^\hdot(\wt{X}/S))
\otimes_{\calo(S)} \Omega^1(S)
$$
obtained by application of the Gauss-Manin connection $\nabla$ to
the relative $2$-form $\Omega$ is called the {\em Kodaira-Spencer
class} of the deformation $\wt{X}/S$.
\end{defn}

Note that the symplectic form $\Omega$ is a cohomology class of the
complex $F^2\Omega^\hdot(\wt{X}/S)$. Since the Gauss-Manin
connection decreases the Hodge filtration at most by $1$, the
Kodaira-Spencer class $\theta$ is well-defined for an arbitrary
symplectic manifold $X$. 

When the symplectic manifold $X$ is admissible, the Kodaira-Spencer
class essentially coincides with the codifferential of the period
map $\Per(\wt{X}):S \to \perdom$. More precisely, the codifferential
\begin{equation}\label{diff.per}
\delta\Per(\wt{X}): \Per^*\Omega^1(\perdom) \to \Omega^1(S)
\end{equation}
of the period map is given by
$$
\delta\Per(\wt{X})(\alpha) = \langle \alpha, \theta \rangle,
$$
where the bundle $\Per^*\Omega^1(\perdom)$ is identified
with the trivial bundle 
$$
\left(H^2(X,F^1\Omega^\hdot(X))\right)^* \otimes \calo_S,
$$ 
$\alpha$ is an arbitrary section of this trivial bundle, and
$\langle \bullet, \bullet \rangle$ means the pairing on the first
factor in
$$
H^2(\wt{X},F^1\Omega^\hdot(\wt{X}/S)) \otimes_{\calo(S)}
\Omega^1(S).
$$
Let now $i:S_0 \hookrightarrow S$ be an elementary extension, let
$\wt{X}_0/S_0$ be a symplectic deformation of a symplectic manifold
$X$, and let
$$
\theta_0 \in H^2(\wt{X}_0,F^1\Omega^\hdot(\wt{X}_0/S_0))
\otimes_{\calo(S_0)} \Omega^1(S_0)
$$
be its Kodaira-Spencer class. Denote by $\eta:i^*\Omega^1(S) \to
\Omega^1(S_0)$ the canonical surjection of the modules of
differentials.

\begin{defn}\label{ks.defn}
By $KS(\theta_0,S)$ we will denote that set of all cohomology
classes
$$
\theta \in H^2(\wt{X}_0,F^1\Omega^\hdot(\wt{X}_0/S_0))
\otimes_{\calo(S_0)} i^*\Omega^1(S)
$$
such that
$$
\eta(\theta) = \theta_0 \in H^2(\wt{X}_0,F^1\Omega^\hdot(\wt{X}_0/S_0))
\otimes_{\calo(S_0)} \Omega^1(S_0).
$$
\end{defn}

Assume that $X$ is admissible, so that for every symplectic
deformation we have the period map, and denote by $p_0:S_0 \to
\perdom$ the period map of the deformation $\wt{X}_0/S_0$. It is the
set $KS(\theta_0,S)$ which we will use as a model for the set
$\Per(p_0,S)$. To do this, notice that every element $p \in
\Per(p_0,S)$ defines an element $\theta(p) \in KS(\theta_0,S)$ by
the formula \eqref{diff.per}. This correspondence is in fact
one-to-one.

\begin{lemma}
The correspondence $p \mapsto \theta(p)$ is a bijection between the
set $\Per(p_0,S)$ and the set $KS(\theta_0,S)$.
\end{lemma}

\proof{} Indeed, both sets are torsors over the group
$$
H^2(\wt{X}_0,F^1\Omega^\hdot(\wt{X}_0/S_0)) \otimes_{\calo(S_0)} I,
$$
where $I \subset \calo(S)$ is the kernel of the map $\calo(S) \to
\calo(S_0)$. For the group $\Per(p_0,S)$, this is obvious from the
exact sequence
$$
\begin{CD}
0 @>>> I @>>> \calo(S) @>>> \calo(S_0) @>>> 0,
\end{CD}
$$
while for the set $KS(\theta_0,S)$ this follows from the exact
sequence \eqref{diff} of the differentials (exact on the left since
the extension $S_0 \subset S$ is elementary).  
\endproof

The advantages of the set $KS(\theta_0,S)$ over the set
$\Per(p_0,S)$ are twofold. Firstly, Definition~\ref{ks.defn} can be
refined so that it takes account of automorphisms of deformation of
$X$. This will be explained in the next section. Secondly, and this
we will state now, Definition~\ref{ks.defn} is essentially local on
$X$. Namely, we have the following obvious fact.

\begin{lemma}
The set $KS(\theta_0,S)$ can be equivalently defined as the set of
all maps
$$
\theta \in \Hom(\calo(S_0)[-2],F^1\Omega^\hdot(\wt{X}_0/S_0)
\otimes_{\calo(S_0)} i^*\Omega^1(S))
$$
in the derived category of sheaves of $\calo(S_0)$-modules on $X$,
such that $\eta(\theta) = \theta_0$. Here $\calo(S_0)$ is considered
as the constant sheaf. \endproof
\end{lemma}

Because of this, we can use Definition~\ref{ks.defn} to restate
Proposition~\ref{rel.main} without the admissibility condition on
the manifold $X$.

\begin{prop}\label{ks.main}
Let $S_0 \subset S$ be an elementary extension of local Artin
algebras.  Let $X$ be a symplectic manifold, let $\wt{X}_0/S_0$ be a
symplectic deformation, and let 
$$
\theta_0 \in \Hom(\calo(S_0)[-2],F^1\Omega^\hdot(\wt{X}_0/S_0)
\otimes_{\calo(S_0)} \Omega^1(S_0))
$$
be its Kodaira-Spencer class. 

Then associating to a symplectic deformation $\wt{X}/S$ its
Kodaira-Spencer class defines a bijection
$$
\Per(\wt{X}_0,S):\Def(\wt{X}_0,S) \cong KS(\theta_0,S).
$$
\end{prop}

\section{Localization, stacks and the proof of the main theorem.}
\label{_Proof_Section_}

\subsection{Groupoids.} We can now begin the proof of
Proposition~\ref{ks.main}, hence also of Theorem~\ref{main}.

Assume given a symplectic manifold $X$, an elementary extension
$i:S_0 \hookrightarrow S$ and a symplectic deformation
$\wt{X}_0/S_0$. Denote by
$$
\theta_0 \in \Hom(\calo(S_0)[-2],F^1\Omega^\hdot(\wt{X}_0/S_0)
\otimes_{\calo(S_0)} \Omega^1(S_0))
$$
the Kodaira-Spencer class of the deformation $\wt{X}_0/S_0$.

To begin with, we will refine the definitions of the sets
$\Def(\wt{X}_0,S)$ and $KS(\theta_0,S)$ and of the map 
$$
\Per(\wt{X}_0,S):\Def(\wt{X}_0,S) \to KS(\theta_0,S)
$$
so as to take into account possible automorphisms of the
deformations $\wt{X} \in \Def(\wt{X}_0,S)$. Such deformations
naturally form a category. Moreover, this category is obviously a
groupoid (i.e. all morphisms are invertible).

\begin{defn}
By $\cDef(\wt{X}_0,S)$ we will denote the groupoid of all
symplectic deformations $\wt{X}/S$ equipped with an isomorphism
$\wt{X} \times_S S_0 \cong \wt{X}_0$.
\end{defn}

The set $\Def(\wt{X}_0,S)$ is by definition the set of objects in
the groupoid $\cDef(\wt{X}_0,S)$.

To construct a groupoid version of the set $KS(\theta_0,S)$, note
that the elements
$$
\theta \in \Hom(\calo(S_0)[-2],F^1\Omega^\hdot(\wt{X}_0/S_0)
\otimes_{\calo(S_0)} \Omega^1(S_0))
$$
naturally classify exact sequences
$$
\begin{CD}
0 @>>> F^1\Omega^\hdot(\wt{X}_0/S_0) \otimes_{\calo(S_0)}
i^*\Omega^1(S) @>>> \bullet @>>> \calo(S_0) @>>> 0
\end{CD}
$$
in the (abelian) category of complexes of sheaves of
$\calo(S_0)$-modules on $X$. Such an element satisfies $\eta(\theta)
= \theta_0$ if and only if there exists a commutative diagram
\begin{equation}\label{ks.diagr}
\begin{CD}
0 @>>> F^1\Omega^\hdot(\wt{X}_0/S_0) \otimes_{\calo(S_0)}
i^*\Omega^1(S) @>>> \bullet @>>> \calo(S_0) @>>> 0\\
@. @V{\id \otimes \eta}VV @VVV @| @.\\
0 @>>> F^1\Omega^\hdot(\wt{X}_0/S_0) \otimes_{\calo(S_0)}
\Omega^1(S_0) @>>> \LL @>>> \calo(S_0) @>>> 0
\end{CD},
\end{equation}
whose the bottom row is the exact sequence corresponding to the
given class 
$$
\theta_0 \in \Hom(\calo(S_0)[-2],F^1\Omega^\hdot(\wt{X}_0/S_0)
\otimes_{\calo(S_0)} \Omega^1(S_0)),
$$
and $\LL$ is the associated extension. 

This motivates the following.

\begin{defn}\label{ks.grp}
By $\KS(\theta_0,S)$ we will denote the groupoid whose objects are
commutative diagrams of the type \eqref{ks.diagr}, and whose
morphisms are maps between these commutative diagrams identical
everywhere except for $\bullet$.
\end{defn}

Again, by definition the set $KS(\theta_0,S)$ is the set of
isomorphism classes of objects in the groupoid
$\KS(\theta_0,S)$. Moreover, the period map
$$
\Per(\wt{X}_0,S):\Def(\wt{X}_0,S) \cong KS(\theta_0,S)
$$
lifts to a functor
\begin{equation}\label{cper}
\cPer(\wt{X}_0,S):\cDef(\wt{X}_0,S) \cong \KS(\theta_0,S).
\end{equation}
This is not immediately obvious, since an element in a first
$\Ext$-group defines a short exact sequence only up to a
non-canonical isomorphism. However, in our case there is a canonical
choice for a short exact sequence \eqref{ks.diagr}. Namely, for
every deformation $\pi:X \to S$ we have a two-step filtration
$\pi^*\Omega^1(S) \subset \Omega^1(X)$ on the sheaf of K\"ahler
differentials $\Omega^1(X)$, with the quotient $\Omega^1(X/S)$. This
filtration induces a filtration on the total de Rham complex
$\Omega^\hdot(X)$. Taking only the top two quotients of this
filtration, we obtain a canonical short exact sequence
\begin{equation}\label{rel.dr}
\begin{CD}
0 @>>> \pi^*\Omega^1(S) \otimes \Omega^\hdot(X/S) @>>> \bullet @>>>
\Omega^{\hdot+1}(X/S) @>>> 0
\end{CD}
\end{equation}
of complexes on $X$ and the corresponding extension class
$$
\eta \in \Ext^1(\Omega^{\hdot+1}(X/S),\pi^*\Omega^1(S) \otimes
\Omega^\hdot(X/S)).
$$
The class $\eta$ essentially induces the Gauss-Manin connection: for
every relative cohomology class 
$$
\alpha \in H^\hdot(X,\Omega^\hdot(X/S)) = \Ext^\hdot(\calo(X),
\Omega^\hdot(X/S)),
$$
the class
$$
\nabla(\alpha) \in H^\hdot(X,\Omega^\hdot(X/S)) \otimes \Omega^1(S)
$$ 
is equal to $\rho\circ\alpha$. Moreover, the sequence \eqref{rel.dr}
is compatible with the Hodge filtration. In particular, we have a
canonical short exact sequence
\begin{equation}\label{rel.dr.bis}
\begin{CD}
0 @>>> \pi^*\Omega^1(S) \otimes F^1\Omega^\hdot(X/S) @>>> \bullet
@>>> F^2\Omega^\hdot(X/S) @>>> 0.
\end{CD}
\end{equation}
Now, by definition of the Kodaira-Spencer class $\theta$ we have
$$
\theta = \nabla(\Omega) = \eta \circ \Omega.
$$ 
Thus the short exact sequence \eqref{ks.diagr} corresponding to the
deformation $X/S$ is obtained by composing the canonical short exact
sequence \eqref{rel.dr.bis} with the map $\calo(X)[2] \to
F^2\Omega^\hdot(X/S)$ given by $\Omega$.

The reader can easily see that this construction is completely
functorial, so that we indeed have a functor \eqref{cper}. We will
call it the {\em period functor}.

Our third and final reformulation of Theorem~\ref{main} is the
following.

\begin{prop}\label{grp.main}
Under the assumption of Proposition~\ref{ks.main}, the period
functor
$$
\cPer(\wt{X}_0,S):\cDef(\wt{X}_0,S) \cong \KS(\theta_0,S)
$$
is an equivalence of categories.
\end{prop}

This proposition implies Proposition~\ref{ks.main}, hence also
Proposition~\ref{rel.main} and Theorem~\ref{main}.

\subsection{Reduction to the affine case.} The following lemma
explains why we introduced the groupoids.

\begin{lemma}
Assume that Proposition~\ref{grp.main} holds for affine symplectic
manifolds $X$. Then it holds for an arbitrary symplectic manifold
$X$.
\end{lemma}

\proof{} For every open subset $U \subset X$, the deformation
$\wt{X}_0/S_0$ induces a symplectic deformation $\wt{U}_0/S_0$ of
the manifold $U$. The collection of groupoids $\cDef(\wt{U}_0,S)$
is a stack on $X$ in Zariski topology.

Moreover, for every open $U \subset X$ any diagram of the type
\eqref{ks.diagr} induces by restriction a diagram of the same type
for the manifold $U$, and the collection of groupoids
$\KS(\theta_0|_U,S)$ is also a stack on $X$ in Zariski topology (it
is to insure this that we have chosen to work with exact sequences
of complexes on $X$ instead of using the derived category).

The period functors $\cPer(\wt{U}_0,S):\cDef(\wt{U}_0,S) \to
\KS(\theta_0|_U,S)$ is a functor between these stack. Since by
assumption it is an equivalence for affine $U \subset X$, it is an
equivalence for an arbitrary $U \subset X$ -- in particular, for $X$
itself.
\endproof

\subsection{The affine case.} 
\label{_affine_ca_Subsection_}
It remains to prove
Proposition~\ref{grp.main} in the case of an affine symplectic
manifold $X$. This is essentially done in the Introduction.
Here we state the same argument in a more refined language. 

Assume given an affine symplectic manifold $X$, an elementary
extension $S_0 \subset S$, and a symplectic deformation
$\wt{X}_0/S_0$ with Kodaira-Spencer class $\theta_0$. 

Since $X$ is affine and smooth, every deformation $\wt{X}/S$ is
trivial as an algebraic manifold. Thus we can choose an isomorphism
$\wt{X}_0 \cong X \times S_0$, and every object $\wt{X} \in
\Ob\cDef(\wt{X}_0,S)$ is isomorphic as an algebraic manifold to $X
\times S$. We introduce the following notation.

\begin{defn}
If a group $G$ acts on a set $N$, then the {\em quotient groupoid}
$N/G$ is the groupoid whose objects are elements $n \in N$ and whose
morphisms are given by
$$
\Hom(n_1,n_2) = \{g \in G\mid g \cdot n_1 = n_2\}, \qquad n_1,n_2
\in N.
$$
\end{defn}

Then the groupoid $\cDef(\wt{X}_0,S)$ is equivalent to the quotient
groupoid
$$
\Sympl/\Aut,
$$
where $\Sympl$ is the set of all relative closed $2$-forms $\Omega
\in \Omega^2(X \times S/S)$ whose restriction to $\wt{X}_0 = X
\times S_0 \subset X \times S$ coincides with the given symplectic
form $\Omega_0 \in \Omega^2(\wt{X}_0/S_0)$, and $\Aut$ is the group
of automorphisms of the algebraic manifold $X \times S$ which
commute with the projection $X \times S \to S$ and which are
identical on $\wt{X}_0 \subset X \times S$.

Every form $\Omega \in \Sympl$ can be represented as
$$
\Omega = \Omega_0 + \beta,
$$
where $\beta \in \Omega^2(X) \otimes I$ is a closed $2$-form on $X$
with values in the ideal $I = \Ker(\calo(S) \to
\calo(S_0)$. Therefore we have a canonical identification $\Sympl =
\Omega^2_{cl}(X) \otimes I$ of $\Sympl$ with the space
$\Omega^2_{cl}X) \otimes I$ of closed $I$-valued $2$-forms on
$X$. Moreover, every automorphism $g \in \Aut$ is an automorphism of
the function ring $\calo(X \times S)$ of the form
$$
g = \id + \xi,
$$
where $\xi \in \T(X) \otimes I$ is an $I$-valued vector field on
$X$. Since $I \subset \calo(S)$ is a square-zero ideal, the group
$\Aut$ is commutative and isomorphic to the abelian group $\T(X)
\otimes I$ -- that is, we have
$$
(\id + \xi_1)\cdot(\id + \xi_2) = \id + \xi_1 + \xi_2
$$
for every $\xi_1,\xi_2 \in \T(X) \otimes I$. Finally, by Cartan
homotopy formula, the action of $\Aut$ on $\Sympl$ is given by
$$
(\id + \xi) \cdot (\Omega_0 + \beta) = \Omega_0 + \beta + d(\Omega_0
\cntrct \xi).
$$
where $d$ is the de Rham differential (all the other terms vanish
since $I^2 = 0$). To sum up, we have
$$
\cDef(\wt{X}_0,S) \cong \Omega^2_{cl}(X) \otimes I/\T(X) \otimes I,
$$
and the action is given by
$$
\xi \cdot \beta = \beta + \Omega_0 \cntrct \xi.
$$
We will now describe the right-hand side of the hypothetical
equivalence \eqref{cper} -- that is, the groupoid $\KS(\theta_0,S)$
-- in a similar way. To do this, notice that since $X$ is affine, we
can replace sheaves on $X$ with the modules of their global
sections. Moreover, the trivialization $\wt{X}_0 \cong X \times S_0$
provides identifications
\begin{align*}
F^1\Omega^\hdot(\wt{X}_0/S_0) \otimes_{\calo(S_0)} i^*\Omega^1(S)
&\cong F^1\Omega^\hdot(X) \otimes_{\C} \Omega^1(S)/I,\\
F^1\Omega^\hdot(\wt{X}_0/S_0) \otimes_{\calo(S_0)} \Omega^1(S_0)
&\cong F^1\Omega^\hdot(X) \otimes_{\C} \Omega^1(S_0),
\end{align*}
of the complexes in the left column of a commutative diagram of type
\eqref{ks.diagr}. 

In every commutative diagram of complexes of this type, both rows
split as exact sequences of graded vector spaces. The only possibly
non-trivial extension data are contained in the differential of the
complex $\bullet$. Denote the set of possible differentials by
$\Diff$. Analogously, every map between two diagrams of the type
\eqref{ks.diagr} must be an automorphism of $\bullet$ which is
upper-triangular with respect to the splitting
$$
\bullet = \left( F^1\Omega^\hdot(\wt{X}_0/S_0) \otimes_{\calo(S_0)}
\Omega^1(S)/I \right) \oplus \calo(S_0).
$$
If we denote the group of all such automorphisms by $\Tr$, then the
groupoid $\KS(\theta_0,S)$ is equivalent to the quotient groupoid
$\Diff/\Tr$. 

To identify the set $\Diff$, note that every possible differential
$\6 \in \Diff$, being a map of $\calo(S_0)$-modules, is completely
determined by its value
$$
\6(1) \in \Omega^2(X) \otimes \Omega^1(S)/I
$$
on the unity $1 \in \calo(S_0)$. Since the differential in the
complex $\LL$ from the bottom row of \eqref{ks.diagr} is fixed,
every two such differentials $\6_1, \6_2 \subset \Diff$ must differ
by a $2$-form
$$
\beta = \6_1(1) - \6_2(1) \in \Omega^2(X) \otimes I \subset
\Omega^2(X) \otimes \Omega^1(S)/I.
$$
Moreover, since every differential $\6 \subset \Diff$ must satisfy
$\6^2=0$, the difference $\beta = \6_1(1) - \6_2(1)$ must be a
closed $I$-valued $2$-form. Thus the set $\Diff$ is a torsor over
the abelian group $\Omega^2_{cl}(X) \otimes I$.

Analogously, every triangular map $g \in \Tr$ is of the form
$$
g = \begin{pmatrix} \id & a\\0 & \id \end{pmatrix}
$$
for some 
$$
a \in \Hom(\calo(S_0),\Omega^1(X) \otimes I) \subset
\Hom(\calo(S_0),\Omega^1(X) \otimes \Omega^1(S)/I),  
$$
and composition of maps $g_1$, $g_2$ simply adds
the associated elements $a_1$, $a_2$. Since $a$ must be a map of
$\calo(S_0)$-modules, it is completely determined by
$$
\alpha = a(1) \in \Omega^1(X) \otimes I.
$$
Therefore $\Tr$ is the abelian group $\Omega^1(X) \otimes I$. Under
these identifications, the action of $\Tr$ on $\Diff$ is given by
$$
\alpha \cdot \beta = \beta + d\alpha,
$$
where $d$ is the de Rham differential.

Having done these identifications, we notice that the period functor
$\cPer:\cDef(\wt{X}_0,S) \to \KS(\theta_0,S)$ is defined on the
level of quotient groupoids by maps
\begin{align*}
\Sympl &\to \Diff\\
\Aut &\to \Tr
\end{align*}
The first map is tautologically an isomorphism $\Omega^2_{cl}(X)
\otimes I \cong \Omega^2_{cl}(X) \otimes I$. The second is the map 
$$
\T(X) \otimes I \to \Omega^1(X) \otimes I
$$
given by $\xi \mapsto \Omega_0 \cntrct \xi$. It is an isomorphism
because the symplectic $2$-form $\Omega_0 \in
\Omega^2(\wt{X}_0/S_0)$ is by assumption non-degenerate.

This finishes the proof of Proposition~\ref{grp.main}, hence, by a
long chains of reductons, also proves Theorem~\ref{main}.
\endproof

\section{Postface}

To finish the paper, we would to give a few comments (perhaps a bit
vague) as to how the theory of symplectic deformations is related to
the usual deformation theory, and what is the relation of this paper
to the known results.

The period map is a purely symplectic phenomenon -- it has no
analogs in the usual deformation theory (although there might be a
useful version for the deformation theory of Calabi-Yau
manifolds). On the other hand, its differential, which we called the
Kodaira-Spencer class, is a fairly general thing. Namely, recall
that for every smooth family $\pi:X \to S$, say over an affine base
$S$, there exist a canonical class
$$
\theta \in \Ext^1(\Omega^1(X/S),\pi^*\Omega^1(S)),
$$
-- the extension class given by the exact sequence
$$
\begin{CD}
0 @>>> \pi^*\Omega^1(S) @>>> \Omega^1(X) @>>> \Omega^1(X/S) @>>> 0
\end{CD}
$$
of differentials for the map $\pi:X \to S$. Since $X/S$ is smooth,
the sheaf $\Omega^1(X/S)$ is flat, and this class can be
reinterpreted as a class in the first cohomology group
$$
H^1(X,\T(X/S) \otimes \pi^*\Omega^1(S)) \cong H^1(X,\T(X/S)) \otimes
\Omega^1(S)
$$
of the relative tangent bundle $\T(X/S)$. It is this class that is
usually called the Kodaira-Spencer class of the deformation $X/S$.

If the deformation $X/S$ is symplectic, then the symplectic form
$\Omega \in \Omega^2(X/S)$ identifies the relative tangent bundle
$\T(X/S)$ with the relative cotangent bundle $\Omega^1(X/S)$ by
means of the correspondence $\xi \mapsto \Omega \cntrct \xi$. This
makes it possible to compare the usual Kodaira-Spencer class
$$
\theta \in H^1(X,\T(X/S)) \otimes \Omega^1(S) \cong
H^1(X,\Omega^1(X/S) \otimes \Omega^1(S))
$$
with the ``symplectic'' Kodaira-Spencer class $\wt{\theta}$
introduced in Definition~\ref{ks.class.defn}. These classes
essentially coincide:

\begin{lemma}
For every symplectic deformation $X/S$, the usual Kodaira-Spencer
class 
$$
\theta = \xi \cntrct \Omega \in H^1(X,\Omega^1(X/S) \otimes
\Omega^1(S))
$$
is obtained from the symplectic Kodaira-Spencer class 
$$
\wt{\theta} \in H^2(X,F^1\Omega^\hdot(X/S)) \otimes \Omega^1(S)) =
H^1(X,F^1\Omega^\hdot(X/S)[1]) \otimes \Omega^1(S))
$$
by the tautological projection
$$
F^1\Omega^\hdot(X/S)[1] \to \Omega^1(X/S).
$$
\end{lemma}

\proof{} Indeed, by definition of the Gauss-Manin connection
$\nabla$ and the usual Kodaira-Spencer class $\xi \in H^1(X,\T(X/S))
\otimes \Omega^1(S)$, for every smooth family $X/S$ and every global
relative $k$-form $\alpha \in \Omega^k(X/S)$ we have
$$
\nabla(\alpha) = \alpha \cntrct \xi.
$$
Applying this to our symplectic deformation $X/S$ and to the
symplectic form $\Omega$ gives the result.
\endproof

In fact, Definition~\ref{ks.defn} of the groupoid $\KS(\wt{X}_0,S)$
is also general and works in the usual deformation theory. Moreover,
the corresponding versions of Proposition~\ref{ks.main} and
Proposition~\ref{grp.main} are also true (and proofs are more or
less the same). However, in the usual case we do not have the period
map, and the groupoid $\KS(\wt{X}_0,S)$ is not easy to describe. In
particular, it might be empty -- that is, there might be a
homological obstruction to the existence of a commutative diagram of
type \eqref{ks.diagr}. This is a very well-known phenomenon. In the
symplectic case, the existence of the period map ensures that (for
admissible manifolds) there are no obstructions for deformation at
any step.

When one tries to describe usual deformations by means of the
associated Kodaira-Spencer class (instead of a period map which does
not exist in general), one enters a closed loop: the class $\theta$,
which theoretically should uniquely define a deformation $X/S$, lies
in the group $H^1(X,\T(X/S))$ which itself depends on the
deformation. The main technical idea in our proof of
Theorem~\ref{main} is to avoid it by going step-by-step through
elementary extensions and using the exact sequence of differentials
$$
\begin{CD}
0 @>>> I @>>> \Omega^1(A)/I @>>> \Omega^1(A/I) @>>> 0
\end{CD}
$$
for an elementary extension $\Spec A/I \subset \Spec A$. This allows
one to describe extension of a deformation $X_0/\Spec (A/I)$ to a
deformation $X/\Spec A$ in terms of the lifting of the
Kodaira-Spencer class from $\Omega^1(A/I)$ to $\Omega^1(A)/I$ -- and
this works, because the module $\Omega^1(A)/I$ is already defined 
over $A/I$. This idea (at least for one-parameter elementary
extensions $\Spec \C[t]/t^k \subset \Spec \C[t]/t^{k+1}$) is due
entirely to Z. Ran \cite{ran}. We believe that it is this technique
that he called {\em the $T_1$-lifting property}.

We would also like to notice that all the obstructions to symplectic
deformations vanish essentially because for an admissible manifold
$X$, the sheaf $R^2\pi_*(\wt{X},F^1\Omega^\hdot(\wt{X}/S))$ is flat
on $S$ for every deformation $\pi:X \to S$. The same thing happens
in the Ran's proof of the Tian-Todorov Lemma -- that is, the lack of
obstructions for deformations of a compact Calabi-Yau manifold is
due to the flatness of some canonically defined sheaves. However,
Ran works with the usual deformations, and he (in the notation as
above) needs the flatness of two sheaves: $\pi_*\Omega^n(\wt{X}/S)$
and $R^1\pi_*\Omega^{n-1}(\wt{X}/S)$ (here $n = \dim X$). This
flatness is provided by Hodge theory. In the symplectic version of
the proof, one would use instead the sheaves
$\pi_*\Omega^2(\wt{X}/S)$ and $R^1\pi_*\Omega^1(\wt{X}/S)$. 

{}From this point of view, the only new thing in our paper is the
following observation: if one agrees to consider deformations that
are {\em a priori} symplectic, then one can combine
$\pi_*\Omega^2(\wt{X}/S)$ and $R^1\pi_*\Omega^1(\wt{X}/S)$ into
$R^2\pi_*F^1\Omega^\hdot(\wt{X}/S)$ -- and the latter sheaf is flat
for a much wider class of manifolds $X$.

Finally, there's another, perhaps more conceptual explanation for
the exceptional role played by the complex $F^1\Omega^\hdot(X)$ in
the symplectic deformation theory. This explanation comes from the
general deformation theory of algebras over an operad, sketched for
example in \cite{gaits1}, \cite{gaits2}. Symplectic manifolds can
not be described by operad. However, a more general class of {\em
Poisson manifolds} admits such a description. The general
deformation theory for algebras over an operad works in a way
completely parallel to the usual one, but the tangent bundle (more
generally, the tangent complex) is replaced by its operadic version.
For the Poisson operad and a symplectic manifold $X$, the Poisson
tangent complex is precisely $F^1\Omega^\hdot(X)$ (independently of
the symplectic form). For a more general Poisson structure, this is
the complex $\Lambda^{\geq 1}\T(X)$ of polyvector fields of degree
$\geq 1$ on $X$, and the differential is given by the commutator
with the Poisson bivector field $\Theta \in \Lambda^2\T(X)$.

\bigskip

{\bf Acknowledgements:} This work is greatly influenced by T.Pantev,
who explained to us the algebraic proof of Bogomolov-Tian-Todorov
theorem.  The second author is grateful to D. Kazhdan,
M. Kontsevich, Yu.I.Manin, S. Merkulov and A. Todorov for
interesting talks on the deformation theory.


\begin{thebibliography}{666}

\bibitem[Ba]{_Barannikov_} 
S. Barannikov,
{\it Generalized periods and mirror symmetry in dimensions $n>3$},
math.AG/9903124, 51 pages.

\bibitem [BK]{_Barannikov_Kontse_} 
S. Barannikov and M. Kontsevich,
{\it Frobenius Manifolds and Formality 
of Lie Algebras of Polyvector Fields},
alg-geom/9710032, 12 pages.


\bibitem[Bea]{_Beauville_} A. Beauville, {\em Varietes
K\"ahleriennes dont la premi\`ere classe de Chern est nulle},
J. Diff. Geom. {\bf 18} (1983), 755--782.

\bibitem[Bo]{bogo} F. Bogomolov, {\em Hamiltonian K\"ahler
manifolds}, Sov. Math. Dokl. {\bf 19} (1978), 1462--1465.

\bibitem[Br]{brisk} E. Brieskorn, {\em Singular elements of
semi-simple algebraic groups}, Actes du Congres International des
Mathematiciens (Nice, 1970), Tome 2, Gauthier-Villars, Paris, 1971,
279--284.

\bibitem[DI]{_Deligne_Illusie_} P. Deligne and L. Illusie, {\em
Relevements modulo $p^2$ et decomposition de complexe de De Rham},
Invent. Math. {\bf 89} (1987), 247--270.

\bibitem[G1]{gaits1} D. Gaitsgory, {\em Operads, Grothendieck
topologies and Deformation theory}, alg-geom/9502010.

\bibitem[G2]{gaits2} D. Gaitsgory, {\em Grothendieck topologies and
deformation theory, II}, alg-geom/9508004.

\bibitem[GPR]{_Grauert_}
Several complex variables. VII. 
Sheaf-theoretical methods in complex analysis. A reprint of {\it Current
problems in mathematics. Fundamental directions. Vol. 74} (Russian), 
Vseross. Inst. Nauchn. i Tekhn. Inform. (VINITI),
Moscow. Edited by H. Grauert, Th. Peternell and R. Remmert. 
Encyclopaedia of Mathematical Sciences, 74. Springer-Verlag,
Berlin, 1994.

\bibitem[KS]{_Kodaira_Spencer_} K. Kodaira and D.C. Spencer, {\em On
deformations of complex structures I, II}, in Kodaira K., Collected
Works, Princeton Univ Press (1975).
 
\bibitem[K]{kawa} Y. Kawamata, {\em Unobstructed deformations -- A
remark on a paper of Z. Ran}, J. Alg. Geom. {\bf 1} (1992),
183--190.

\bibitem[R]{ran} 
Z. Ran, {\em Deformations of manifolds with torsion
or negative canonical bundles}, J. Alg. Geom. {\bf 1} (1992),
279--291. 

\bibitem[T]{tian} G. Tian, {\em Smoothness of the universal
deformation space of compact Calabi-Yau manifolds and its
Petersson-Weil metric}, in {\em Math. Aspects of String Theory},
S.-T. Yau, ed., Worlds Scientific, 1987, 629--646.

\bibitem[To]{todo} A. Todorov, {\em The Weil-Petersson geometry of
the moduli space of $SU(n \geq 3)$ (Calabi-Yau) manifolds}, Comm.,
Math. Phys. {\bf 126} (1989), 325--346.

\end{thebibliography}
\end{document}